\documentclass{ifacconf}

\usepackage{graphicx}      
\usepackage{natbib}        
\usepackage{xcolor}
\usepackage{amssymb}
\usepackage{mathtools}
\usepackage{tikz-cd}
\usepackage{textcomp}

\newcommand{\Cont}{\mathcal{C}}
\newcommand{\CL}{\Cont_{\textnormal{loc}}}

\newcommand{\CLKA}{\Cont_{\textnormal{loc}}^{k,\alpha}}
\newcommand{\T}{\mathsf{T}}
\newcommand{\D}{\mathsf{D}}
\newcommand{\GL}{\mathsf{GL}}

\newcommand{\Stab}{\mathcal{S}_n}
\newcommand{\Stabp}{\mathcal{S}_{n-1}}
\newcommand{\Stabc}{\bar{\mathcal{S}}_n}

\newcommand{\Nonr}{\mathcal{N}_n}
\newcommand{\Nonp}{\mathcal{N}_{n-1}}
\newcommand{\Goodf}{\mathcal{G}_{\textnormal{fix}}}
\newcommand{\Goodp}{\mathcal{G}_\textnormal{per}}

\newcommand{\Xfinf}{\mathfrak{X}_{\textnormal{fix}}^\infty}
\newcommand{\Xpinf}{\mathfrak{X}_{\textnormal{per}}^\infty}

\newcommand{\N}{\mathbb{N}}
\newcommand{\Z}{\mathbb{Z}}
\newcommand{\R}{\mathbb{R}}
\newcommand{\Grp}{\mathbb{R}}
\newcommand{\C}{\mathbb{C}}
\newcommand{\slot}{\,\cdot\,} 

\newcommand{\concept}[1]{\textit{#1}}

\newcommand{\david}[1]{[{\color{magenta}(D.H.) \color{cyan}\textsc{#1}}]}

\begin{document}
\begin{frontmatter}

\title{Generic Properties of Koopman Eigenfunctions for Stable Fixed Points and Periodic Orbits\thanksref{footnoteinfo}}

\thanks[footnoteinfo]{Kvalheim and Revzen were supported by ARO award W911NF-14-1-0573 to Revzen and by the ARO under the Multidisciplinary University Research Initiatives (MURI) Program, award W911NF-17-1-0306 to Revzen.
Kvalheim was also supported by the ARO under the SLICE
MURI Program, award W911NF-18-1-0327.
Hong was supported in part by the Dean's Fund for Postdoctoral Research
of the Wharton School.}

\author[First]{Matthew D. Kvalheim}
\author[Second]{David Hong}
\author[Third]{Shai Revzen}

\address[First]{Department of Electrical and Systems Engineering, University of Pennsylvania,
   Philadelphia, PA 19104 USA (e-mail: kvalheim@seas.upenn.edu).}
\address[Second]{Department of Statistics, University of Pennsylvania,
   Philadelphia, PA 19104 USA (e-mail: dahong67@wharton.upenn.edu).}
\address[Third]{Department of Electrical Engineering and Computer Science, Ecology and Evolutionary Biology Department, Robotics Institute, University of Michigan,
   Ann Arbor, MI 48109 (e-mail: shrevzen@umich.edu)}

\begin{abstract}                
Our recent work established existence and uniqueness results for $\Cont^k$ (actually $\CLKA$) linearizing semiconjugacies for $\Cont^1$ flows defined on the entire basin of an attracting hyperbolic fixed point or periodic orbit \citep{kvalheim2019existence}.
Applications include (i) improvements, such as uniqueness statements, for the Sternberg linearization and Floquet normal form theorems, and (ii) results concerning the
existence, uniqueness, classification, and convergence of various quantities appearing in the ``applied Koopmanism'' literature, such as principal eigenfunctions, isostables, and Laplace averages.

In this work we consider the broadness of applicability of these results with an emphasis on the Koopmanism applications.
In particular we show that, for the flows of ``typical'' $\Cont^\infty$ vector fields having an attracting hyperbolic fixed point or periodic orbit with a fixed basin of attraction, the $\Cont^\infty$ Koopman eigenfunctions can be completely classified, generalizing a result known for analytic eigenfunctions of analytic systems.
\end{abstract}


\end{frontmatter}

\section{Introduction}
Linear dynamical systems and control systems are very well understood, in contrast with their nonlinear counterparts.
Most models of real-world systems are, unfortunately, nonlinear.
Thus, any means for applying linear systems techniques to the analysis and synthesis of nonlinear systems is of general interest in both scientific and engineering applications.

A common approach is to \emph{approximate} a nonlinear system as a linear system near some nominal trajectory and apply linear systems techniques to the approximation \citep[Sec.~4.3,~12.2]{khalil2002nonlinear}.
While this approach works well in many situations, it is inherently \emph{local} and often fails if the system is sufficiently far from the nominal trajectory.
A recent alternative approach seeks linear representations of nonlinear systems that are instead \emph{global} and \emph{exact}.
This is the approach taken in the applied Koopman operator theory literature, initiated largely by \citet{mezic1994geometrical, mezic2004comparison, mezic2005spectral}, around 70 years after Koopman's seminal work \citep{koopman1931hamiltonian}.

The Koopman operator of a (nonlinear) dynamical system is an infinite-dimensional linear operator that acts on scalar-valued functions of state, or \concept{observables}, by evolving them via the underlying dynamics.
Since this operator is linear, one can discuss its spectral theory, and its spectral objects often have dynamical relevance.
In particular, Koopman eigenfunctions are observables that evolve linearly under the dynamics; the dynamics and control of observables spanned by Koopman eigenfunctions is thus governed by linear systems theory.
We emphasize that such a reduction is both \emph{exact} and \emph{global}; furthermore, given enough independent eigenfunctions one obtains an exact, global change of coordinates transforming the nonlinear system into a linear one.

Thus, methods to identify Koopman eigenfunctions are of interest, and the numerical computation of such eigenfunctions is an active research area.\footnote{
    Due to space constraints we  mention only the review \citet{budivsic2012applied} here; see the references in \cite{kvalheim2019existence} for many additional examples of this literature.
}
The body of work most relevant to the present paper concerns the numerical computation of  \concept{isostables} \citep{mauroy2013isostables} and \concept{isostable coordinates} \citep{wilson2016isostable, shirasaka2017phase, wilson2018greater, monga2019phase} for dynamical systems having an asymptotically stable equilibrium  or limit cycle; these objects can be expressed in terms of Koopman eigenfunctions as discussed in \citet{kvalheim2019existence}.
Isostables and isostable coordinates are useful tools for nonlinear model reduction, and it has been proposed that these objects could prove useful in real-world applications such as treatment design for Parkinson's disease, migraines, cardiac arrhythmias \citep{wilson2016isostable-pde}, and jet lag \citep{wilson2014energy}.

In analyzing the theoretical properties of any algorithm for computing some quantity, it is desirable to know whether the computation is well-posed \citep{hadamard1902problemes}, and in particular whether the quantity in question \emph{exists} and is \emph{uniquely determined}.
An existence and uniqueness theory for isostables, isostable coordinates, and more general Koopman eigenfunctions is thus desirable.
In the context of attracting equilibria and limit cycles, some existence results can be obtained by invoking Hartman-Grobman type linearization theorems \citep{lan2013linearization,kvalheim2018global}.
On the other hand, it seems that uniqueness was less well understood, with an exception for the case of analytic eigenfunctions for analytic dynamical systems having a nonresonant linearization \citep{mauroy2013isostables, mezic2019spectrum}.

Our recent work \citep{kvalheim2019existence} filled much of the gap by establishing existence and uniqueness results for $\CLKA$ linearizing semiconjugacies for $\Cont^1$ dynamical systems, of which Koopman eigenfunctions are a special case.
In particular, we obtained uniqueness results for $\Cont^k$ Koopman eigenfunctions; we also obtained $\Cont^k$ existence results that, to the best of our knowledge, are stronger than those appearing elsewhere in the literature for $2\leq k \leq \infty$.
We obtained a particularly strong result for the case $k = \infty$: the $\Cont^\infty$ eigenfunctions admit a complete classification for $\Cont^\infty$ dynamical systems satisfying a nondegeneracy condition.
The conclusion of this classification result yields much information about the eigenfunctions, so one would naturally like to understand how often its hypotheses hold.

The contribution of the present work is to show that the classification, existence, and uniqueness results for $\Cont^\infty$ eigenfunctions in \cite{kvalheim2019existence} in fact hold for ``typical'' $\Cont^\infty$ vector fields having an asymptotically stable equilibrium or periodic orbit with a fixed basin of attraction, where ``typical'' means for sets of $\Cont^\infty$ vector fields which are open and dense in suitable topologies.

The remainder of the paper is organized as follows.
In \S \ref{sec:generic} we prove our main result, Theorem \ref{th:main-thm}, after some preliminary definitions and lemmas.
In \S\ref{sec:old-results} we discuss the implications of Theorem \ref{th:main-thm} for the results of \cite{kvalheim2019existence} relevant to Koopman eigenfunctions.
Finally, Appendix \ref{app:polynomial} contains background on symmetric polynomials for the convenience of the reader.

\section{Main results}\label{sec:generic}

This section contains our main result, Theorem  \ref{th:main-thm}.
But first, we need some preliminary definitions and lemmas.

The following definition is \citet[Def.~1]{kvalheim2019existence} and is essentially an asymmetric version of definitions appearing in \cite{sternberg1957local, sell1985smooth}.
By some abuse of notation we also apply this definition to real matrices by viewing them as complex matrices with real entries; when discussing eigenvalues and eigenvectors of a linear self-map or matrix in this work, we always mean eigenvalues and eigenvectors of its complexification. 
By a further abuse of notation we also apply this definition to (the complexifications of) general linear self-maps of finite-dimensional vector spaces.

\begin{defn}[$(X,Y)$ $k$-nonresonance]\label{def:nonres}
    Let $X \in \C^{d\times d}$ and $Y \in \C^{n\times n}$ be matrices with eigenvalues $\mu_1,\ldots,\mu_d$ and $\lambda_1,\ldots,\lambda_n$, respectively, repeated with multiplicities.
    For any $k \in \N_{\geq 1}\cup \{\infty\}$, we say that $(X,Y)$ is \concept{$k$-nonresonant} if, for any $i\in \{1,\ldots, d\}$ and any $m=(m_1,\ldots, m_n) \in \N^n_{\geq 0}$ satisfying $2\leq m_1 + \cdots + m_n < k+1$,
    \begin{equation}\label{eq:nonres-lem}
	    \mu_i \neq \lambda_1^{m_1}\cdots \lambda_n^{m_n}.
    \end{equation}
    (Note that this condition vacuously holds if $k = 1$; i.e., any two matrices are $1$-nonresonant.)
\end{defn}

For $n\in \N_{\geq 1}$ let $\Nonr\subset \R^{n\times n}$ be the set of $n\times n$ real matrices $A$ \emph{with distinct eigenvalues} such that $(A,A)$ is $\infty$-nonresonant; by abuse of notation we also apply this definition to linear self-maps of general $n$-dimensional real vector spaces.
Denoting by $\GL(n,\R)\subset \R^{n\times n}$ the invertible matrices, it follows from Def.~\ref{def:nonres} that $\Nonr \subset \GL(n,\R)$ (since $0 = 0^m$ for all $m\in \N)$.
Below we use the notation $\exp\colon \R^{n\times n}\to \GL(n,\R)$, $\exp(A)\coloneqq e^A$, when convenient.

\begin{lem}\label{lem:nonres-full-measure}
	$\R^{n\times n}\setminus \Nonr$ and $\R^{n\times n}\setminus \exp^{-1}(\Nonr)$ both have Lebesgue measure zero.
\end{lem}

\begin{pf}
	From Def.~\ref{def:nonres}, matrices in $\R^{n\times n} \setminus \Nonr$ have eigenvalues $\lambda_1,\dots,\lambda_n$	that satisfy:	(i) $\lambda_j = \lambda_k$ for some $j \neq k$ or (ii) $\lambda_i = \lambda_1^{m_1} \cdots \lambda_n^{m_n}$	for some $i \in \{ 1, \dots, n \}$ and $(m_1,\dots,m_n) \in \mathcal{M}_n$ where
	\begin{equation*}
	    \mathcal{M}_n \coloneqq \{(m_1,\dots,m_n) \in \N^n_{\geq 0} : m_1 + \cdots + m_n \geq 2\}.
	\end{equation*}
	Condition (i) is equivalent to
	\begin{equation*}
	    0 = f(\lambda_1,\dots,\lambda_n) \coloneqq \prod_{j \neq k} (\lambda_j - \lambda_k),
	\end{equation*}
	and condition (ii) is equivalent to
	\begin{equation*}
	    \exists m \in \mathcal{M}_n : 0 = g_m(\lambda_1,\dots,\lambda_n),
	\end{equation*}
	where
	\begin{equation*}
	    g_m(\lambda_1,\dots,\lambda_n)  \coloneqq \prod_i
	\prod_{\sigma \in S_n} (\lambda_i - \lambda_1^{m_{\sigma(1)}} \cdots \lambda_n^{m_{\sigma(n)}})
	\end{equation*}
	and $S_n$ is the group of permutations $\sigma$ of $\{ 1,\dots,n \}$.
	Since $f$ and $g_m$ (for any $m$) are symmetric polynomials in the eigenvalues, they are also expressible as polynomials $F, G_m: \R^{n\times n} \to \R$ in the matrix entries.
	This follows from the fundamental theorem of symmetric polynomials	and Vieta's theorem	by recalling that the eigenvalues are roots of the characteristic polynomial		whose coefficients are polynomials in the matrix entries (for more details see Appendix \ref{app:polynomial}).
	None of them are identically zero since, e.g.,
	\begin{equation*}
	    F(\operatorname{diag}(\lambda_1,\dots,\lambda_n)) \equiv f(\lambda_1,\dots,\lambda_n) \not\equiv 0,
	\end{equation*}
	and likewise for each $G_m$.
	As a result,
	\begin{equation*}
	   \R^{n\times n} \setminus \Nonr = F^{-1}(0) \cup \bigcup_{m \in \mathcal{M}_n} G_m^{-1}(0)
	\end{equation*}
	is a countable union of measure zero sets and so is also measure zero.
	In more detail: each set in the union is measure zero since (i) polynomials are real analytic functions and (ii) the zero set of a real analytic function which is not identically zero has measure zero \citep{mityagin2015zero}.
	Defining the real analytic functions $\tilde{F}\coloneqq F\circ \exp$ and $\tilde{G}_m\coloneqq G_m\circ \exp$, 
	\begin{equation*}
   \R^{n\times n}\setminus \exp^{-1}(\Nonr) = \tilde{F}^{-1}(0) \cup \bigcup_{m \in \mathcal{M}_n} \tilde{G}_m^{-1}(0)
	\end{equation*}
	is measure zero by the same reasoning.
	 $\square$
\end{pf}

Let $\Stab\subset \Stabc\subset \R^{n\times n}$ denote the sets of $n\times n$ real matrices whose eigenvalues belong to the open and closed unit disks in $\C$, respectively.
Given any matrix $X\in \C^{n\times n}$, we define its \concept{spectral radius} $\rho(X)\coloneqq \max_{\mu\in \textnormal{spec}(X)}|\mu|$, where $\textnormal{spec}(X)\subset \C$ denotes the set of eigenvalues of a matrix $X$.
By abuse of notation we also apply the preceding two definitions to linear self-maps of general $n$-dimensional real vector spaces; there is no ambiguity since eigenvalues do not depend on a choice of basis.

\begin{lem}\label{lem:open-matrices}
	$\Stab \cap \Nonr$ is open in $\Stab$ and $\Stab$ is open in $\R^{n\times n}$.
\end{lem}

\begin{pf}
	Fix any $A\in \Stab\cap \Nonr \subset \GL(n,\R)$.
	Since $A\in \Stab$, there exists $k\in \N_{\geq 2}$ such that $\rho(A^{-1})\rho(A)^{k+1}<1$.
	It follows from Def.~\ref{def:nonres} that $\infty$-nonresonance of $(B,B)$ is implied by (i) $k$-nonresonance of $(B,B)$, which implies that $B$ is invertible, and (ii) $\rho(B^{-1})\rho(B)^{k+1}<1$.
	Since the inverse and eigenvalues of a matrix depend continuously on the matrix \citep[p.~53]{palis1982geometric}, the set of matrices satisfying each of these two conditions is open in $\Stab$.
	Similarly, the set of matrices having distinct eigenvalues is also open in $\Stab$.
	Hence $A$ has a neighborhood in $\Stab$ contained in $\Nonr$;
	since $A$ was arbitrary, $\Stab\cap \Nonr$ is open in $\Stab$.
	Continuity of eigenvalues also directly implies openness of $\Stab$ in $\R^{n\times n}$.
	$\square$
\end{pf}

\begin{lem}\label{lem:dense-matrices}
$\exp^{-1}(\Stab\cap \Nonr)$	 is dense in $\exp^{-1}(\Stabc)\subset \R^{n\times n}$.
\end{lem}
\begin{pf}
Note that $\exp^{-1}(\Stab) \subset \exp^{-1}(\Stabc)$ are the sets of matrices having only eigenvalues with negative and nonpositive real parts, respectively.
Examination of the real canonical form of matrices $A\in \exp^{-1}(\Stabc)$ reveals that $\exp^{-1}(\Stab)$ is dense in $\exp^{-1}(\Stabc)$.
Lem.~\ref{lem:open-matrices} and continuity of $\exp$ imply that $\exp^{-1}(\Stab)$ is open in $\R^{n\times n}$, so Lem.~\ref{lem:nonres-full-measure} implies that $\exp^{-1}(\Stab\cap \Nonr)$ is dense in $\exp^{-1}(\Stab)$ and thus (by the preceding sentence) also in $\exp^{-1}(\Stabc)$. $\square$ 
\end{pf}

Recall that a $\Cont^k$ ($k\in \N_{\geq 1}$) \concept{flow} on a smooth manifold $Q$ is a $\Cont^k$ map $\Phi\colon Q\times \Grp \to Q$ satisfying $\Phi^0 = \textnormal{id}_Q$ and $\Phi^{t+s}=\Phi^t\circ \Phi^s$ for all $t,s\in \Grp$, where $\Phi^t\coloneqq \Phi(\slot,t)$.
(A $\Cont^k$ map is one which has continuous mixed partial derivatives up to order $k$ in local coordinates.)
As a typical example, an ordinary differential equation (ODE)
\begin{equation}\label{eq:ode}
\frac{dx}{dt} = f(x)
\end{equation}
defined by a $\Cont^k$ vector field which is \concept{complete} \citep[p.~215]{lee2013smooth} generates a unique $\Cont^k$ flow $\Phi$, where $t\mapsto \Phi^t(x_0)$ is the unique solution to \eqref{eq:ode} with initial condition $\Phi^0(x_0)=x_0$.
(Any $\Cont^1$ vector field is complete when restricted to the basin of attraction of a compact asymptotically stable set.)

We use the following notation in the remainder of this paper.
Given a differentiable map $F\colon M\to N$ between smooth manifolds, $\D_x F$ denotes the derivative of $F$ at the point $x\in M$.
(Recall that $\D_x F\colon \T_x M \to \T_{F(x)}N$ is a linear map between tangent spaces \citep{lee2013smooth}, which can be identified with the Jacobian of $F$ evaluated at $x$ in local coordinates.)
In particular, given a $\Cont^1$ flow $\Phi\colon Q\times \Grp \to Q$ and fixed $t \in \Grp$, we write $\D_{x}\Phi^t\colon \T_{x}Q\to \T_{\Phi^t(x)}Q$ for the derivative of the time-$t$ map $\Phi^t\colon Q\to Q$ at the point $x$.

We need some additional notation for our main result.
Let $Q$ be a smooth $n$-dimensional manifold with $n\geq 1$. 
Let $\Xfinf(Q)$ and $\Xpinf(Q)$ be the sets of $\Cont^\infty$ vector fields $f$ whose flows possess an asymptotically stable fixed point $x_f$ with basin $Q$ and asymptotically stable nonstationary periodic orbit $\Gamma_f$ with basin $Q$, respectively.
For the case of nonstationary periodic orbits we assume that $\dim(Q)\geq 2$.
We use the notation $\Phi_f$ for the flow of such a vector field $f$.
Given $f\in \Xpinf$, we let $x_f\in \Gamma_f$ be an arbitrary point and $\tau_f>0$ be the period of $\Gamma_f$; if $\Gamma_f$ is hyperbolic, we let $E^s_{x_f}$ be the unique $\D_{x_f}\Phi_f^{\tau_f}$-invariant complement to $\textnormal{span}(\{f(x_f)\})$.
Let $\Goodf \subset \Xfinf$ and $\Goodp \subset \Xpinf$ denote the ``good'' vector fields such that every $f\in \Goodf$  satisfies $e^{\D_{x_f}f}= \D_{x_f}\Phi_{f}^1 \in \Stab \cap \Nonr$ and such that the periodic orbit for each $g\in \Goodp$ is hyperbolic and satisfies $\D_{x_g}\Phi_{g}^{\tau_g}|_{E^s_{x_g}} \in \Stabp \cap \Nonp$.

The theorem below is our main result.
We refer the reader to \citet[Ch.~2]{hirsch1976differential} for the definitions of the $\Cont^k$ Whitney (strong) and compact-open (weak) topologies, but the theorem's effective meaning is clear from its proof.
\begin{thm}\label{th:main-thm}
	$\Goodf$ (resp. $\Goodp$) is open in $\Xfinf(Q)$ (resp. $\Xpinf(Q)$) with respect to the $\Cont^1$ compact-open topology and dense in $\Xfinf(Q)$ (resp. $\Xpinf(Q)$) with respect to the $\Cont^\infty$ Whitney topology.
\end{thm}

\begin{rem}
Many results proved in \cite{kvalheim2019existence}, including those recapitulated in the following \S \ref{sec:old-results}, hold for flows of $\Cont^\infty$ vector fields belonging to $\Goodf$ or $\Goodp$.
Thus, a ``typical'' $\Cont^\infty$  vector field in $\Xfinf(Q)$ or $\Xpinf(Q)$ satisfies the hypotheses of those results.
\end{rem}
\begin{rem}\label{rem:meas}
Despite the suggestive statement of Lem.~\ref{lem:nonres-full-measure}, we have not attempted to formalize ``typical'' in a measure-theoretic sense in Theorem \ref{th:main-thm} due to the apparent lack of natural definitions of ``measure zero'' subsets of $\Xfinf$ and $\Xpinf$.
In this direction, it would be interesting to know whether ``typical'' could be interpreted in a stronger sense using the framework of \concept{prevalence} \citep{ott2005prevalence}.
\end{rem}

\begin{pf}
	We prove the theorem for $\Xfinf$; the case of $\Xpinf$ is handled similarly using Floquet theory.
	The $\Xfinf$ statements hold vacuously if  $\Xfinf(Q) = \varnothing$; if $\Xfinf(Q) \neq \varnothing$ then $Q$ is diffeomorphic to $\R^n$ \citep{wilson1967structure}, so we may henceforth assume that $Q = \R^n$ and that $x_f = 0$.

	\emph{Density} --- %
	Let $f\in \Xfinf(\R^n)$ be arbitrary and let $U\subset Q$ be a precompact open neighborhood of $0\,(=x_f)$.
	Let $\varphi\colon \R^n\to [0,\infty)$ be a $\Cont^\infty$ function equal to $1$ on a neighborhood of $0$ and having support contained in $U$.
	Since $0$ is asymptotically stable, $\D_0 f\in \exp^{-1}(\Stabc)$.
	Lem.~\ref{lem:dense-matrices} implies the existence of a sequence $(A_n)_{n\in \N}$ of matrices with $A_n\to \D_0 f$ and with $e^{A_n}\in \Stab \cap \Nonr$ for all $n$.
	We now define a sequence $(g_n)_{n\in \N}$ of $\Cont^\infty$ vector fields with $\D_{0}\Phi_{g_n}^1 = e^{A_n} \in \Stab \cap \Nonr$ via $$g_n(x)\coloneqq f(x) + \varphi(x) (A_n-\D_{0}f)\cdot x.$$
	All derivatives of the $g_n$ converge uniformly to those of $f$ on $U$, and $g_n$ is equal to $f$ on $\R^n\setminus U$, so $g_n$ converges to $f$ in the $\Cont^\infty$ Whitney topology.
	Note that $g_n(0) = f(0) = 0$ for all $n$.
	It remains only to prove that $g_n\in \Xfinf(\R^n)$ for all $n$ sufficiently large, i.e., that $0$ is globally asymptotically stable for $g_n$ for large $n$; this follows from a general result of \citet[Thm~2.2]{smith1999perturbation}.

	\emph{Openness} --- %
	Fix any vector field $f\in \Goodf\subset \Xfinf(\R^n)$, so that $\D_{0}\Phi_f^1\in \Stab \cap \Nonr$.
	Let $(g_n)_{n\in \N}$ be a sequence of vector fields in $\Xfinf(\R^n)$ converging to $f$ in the $\Cont^1$ compact-open topology; i.e., $g_n$ and $\D g_n$ converge to $f$ and $\D f$ uniformly on compact sets.
	Since the $\Cont^1$ compact-open topology can be given the structure of a Banach space, the (Banach space version of the) implicit function theorem implies that $x_{g_n}\to 0$ and hence $\D_{x_{g_n}}\Phi_{g_n}^1 = e^{\D_{x_{g_n}}g_n}\to e^{\D_{0}f} = \D_{0}\Phi_{f}^1$.
	It follows from Lem.~\ref{lem:open-matrices} that $\D_{x_{g_n}}\Phi_{g_n}^1\in \Stab \cap \Nonr$ and hence $g_n\in \Goodf$ for all $n$ sufficiently large.
	Since $\Xfinf$ with the $\Cont^1$ compact-open topology is first countable, this implies the desired openness statement and completes the proof. $\square$
\end{pf}

\section{Implications for the Existence and Uniqueness of Koopman Eigenfunctions} \label{sec:old-results}
The remainder of this paper describes the implications of Theorem \ref{th:main-thm} for the results of \citet[Sec.~3.2--3.3]{kvalheim2019existence}.
Those results were stated in terms of $\CLKA$ functions; here we discuss only the simpler case $\Cont^k = \CL^{k,0}$.

\subsection{Koopman eigenfunctions}\label{sec:p-eigs}

Given a $\Cont^1$ flow $\Phi\colon Q\times \Grp \to Q$, where $Q$ is a smooth manifold, we say that $\psi \colon Q\to \C$ is a \concept{Koopman eigenfunction} with \concept{eigenvalue} $\mu\in \C$ if $\psi$ is not identically zero and
\begin{equation}\label{eq:koopman-efunc}
\forall t \in \Grp\colon \psi \circ \Phi^t = e^{\mu t} \psi.
\end{equation}
The following generalizes the definitions for linear systems given in \citet[Def.~2.2--2.3]{mohr2016koopman}.

	\begin{defn}\label{def:principal-eigenfunction}
		If $Q$ is the basin of an asymptotically stable fixed point $x_0\in Q$ for $\Phi$, we say that an eigenfunction $\psi \in \Cont^1(Q,\C)$ is a \concept{principal} eigenfunction if $\psi(x_0) = 0$ and $\D_{x_0}\psi \neq 0$.
		If instead $Q$ is the basin of an asymptotically stable periodic orbit with image $\Gamma \subset Q$ for $\Phi$, we say that an eigenfunction $\psi \in \Cont^1(Q,\C)$ is a \concept{principal} eigenfunction if $\psi(x_0)= 0$ and $\D_{x_0}\psi \neq 0$ for all $x_0\in \Gamma$.
	\end{defn}

\subsection{Principal eigenfunctions for fixed points and periodic orbits}\label{sec:p-eig-results}

Given a (real or complex) linear self-map $Y\colon V\to V$, we say that a linear map $w\colon V\to \C$ is a \concept{left eigenvector} of $Y$ with \concept{eigenvalue} $\lambda\in \C$ if $wY = \lambda w$. 
Differentiating \eqref{eq:koopman-efunc} and using the chain rule immediately yields Prop.~ \ref{prop:p-eig-evec-point} and \ref{prop:p-eig-evec-cycle}, which have appeared in the literature (see, e.g., the proof of \cite[Prop.~2]{mauroy2016global}).
\begin{prop}\label{prop:p-eig-evec-point}
	Let $x_0$ be an asymptotically stable fixed point of the flow of a $\Cont^1$ vector field $f$ with basin $Q$. If $\psi\in \Cont^1(Q,\C)$ is a principal Koopman eigenfunction for the flow of $f$ with eigenvalue $\mu\in \C$, then $\D_{x_0}\psi$ is a left eigenvector of $\D_{x_0}f$ with eigenvalue $\mu$.
\end{prop}

\begin{prop}\label{prop:p-eig-evec-cycle}
	Let $\Gamma$ be the image of an asymptotically stable  $\tau$-periodic orbit of the $\Cont^1$ flow $\Phi$ with basin $Q$.\footnote{We always mean that $\tau$ is the \emph{minimal} period in ``$\tau$-periodic orbit''.}
	If $\psi\in \Cont^1(Q,\C)$ is a principal Koopman eigenfunction for $\Phi$ with eigenvalue $\mu\in \C$, then for any $x_0 \in \Gamma$, $\D_{x_0}\psi$ is a left eigenvector of $\D_{x_0}\Phi^\tau$ with eigenvalue $e^{\mu \tau}$.
\end{prop}

The following result follows from \citet[Rem.~3, Ex.~1, Prop.~6]{kvalheim2019existence}.
The condition ``$\D_{x_0}^i R = 0$ for all $0\leq i < k$'' should be interpreted to mean that, in local coordinates, $R$ and all of its mixed partial derivatives of order less than $k$ vanish at $x_0$. 
This does not depend on the choice of local coordinates; see \cite[Sec.~1.3.3]{kvalheim2019existence}.
The $\Cont^k$ compact-open (weak) topology \citep[Ch.~2]{hirsch1976differential} on functions referred to below is the topology of $\Cont^k$-uniform convergence on compact subsets.

\begin{prop}\label{prop:koopman-cka-fix}
	Let $f$ be a $\Cont^1$ vector field on $Q$ with $Q$ the basin of an attracting hyperbolic equilibrium $x_0\in Q$ for the flow of $f$, where $n\coloneqq \dim(Q)\geq 1$.
	Fix $k \in \N_{\geq 1}\cup\{\infty\}$ and assume the spectral radius $\rho\left(e^{\D_{x_0}f}\right)<1$ satisfies 
  	$$|e^\mu|>\left(\rho\left(e^{\D_{x_0}f}\right)\right)^{k}$$
	in all of the following statements (with $\left(\rho\left(e^{\D_{x_0}f}\right)\right)^{\infty}\coloneqq 0$).
	
	\emph{\textbf{Uniqueness of Koopman eigenvalues and principal eigenfunctions.}}
	Let $\psi_1\in \Cont^k(Q,\C)$ be any Koopman eigenfunction with eigenvalue $\mu$.
	\begin{enumerate}
		\item\label{item:uniq-thm-1} Then there exists $m = (m_1,\ldots,m_n)\in \N_{\geq 0}^n$ such that $$\mu =  m\cdot \lambda,$$ where $\lambda_1,\ldots,\lambda_n$ are the eigenvalues of $\D_{x_0}f$ repeated with multiplicities and $\lambda\coloneqq (\lambda_1,\ldots,\lambda_n)$.
		\item\label{item:uniq-thm-2} Assume that $\psi_1$ is a principal eigenfunction so that $\mu \in \textnormal{spec}(\D_{x_0}f)$, and assume that $\mu \neq m\cdot \lambda$ for any $m\in \N^n_{\geq 0}$ with $2\leq \sum_i m_i\leq k$. 
		Then $\psi_1$ is uniquely determined by $\D_{x_0}\psi_1$, and if $\mu$ and $\D_{x_0}\psi_1$ are real, then $\psi_1\colon Q\to \R\subset \C$ is real.
		In particular, if $\mu$ is an algebraically simple eigenvalue of (the complexification of) $\D_{x_0}f$ and if $\psi_2$ is any other principal eigenfunction with eigenvalue $\mu$, then there exists $c\in \C\setminus \{0\}$ such that $$\psi_1 = c\psi_2.$$
	\end{enumerate}

	\emph{\textbf{Existence of principal eigenfunctions.}} Assume $f \in \Cont^k$ and that $\mu \neq m\cdot \lambda$ for any $m\in \N^n_{\geq 0}$ with $2\leq \sum_i m_i\leq k$. 
	Let $w$ be a left eigenvector of $\D_{x_0}f$ with eigenvalue $\mu$.
	\begin{enumerate}
		\item Then there exists a unique principal eigenfunction $\psi \in \Cont^k(Q,\C)$  with eigenvalue $\mu$ satisfying  $\D_{x_0}\psi = w$.
		\item  In fact, if $\Phi$ is the flow of $f$ and $P\in \Cont^k(Q,\C)$ is any ``approximate eigenfunction'' satisfying
		\begin{equation}\label{eq:Koop-main-approx}
		P \circ \Phi^1 = e^\mu P + R
		\end{equation}
		with $\D_{x_0}^i R = 0$ for all $0\leq i < k$ and $\D_{x_0}P = w$, then
		\begin{equation}\label{eq:Koop-main-converge}
		\psi = \lim_{t\to \infty} e^{-\mu t} P \circ \Phi^t
		\end{equation}
		with convergence in the $\Cont^k$ compact-open topology.
	\end{enumerate}
\end{prop}

\begin{rem}[the $\Cont^\infty$ case]\label{rem:point-cinf-case}
	In the case $k = \infty$, the \concept{spectral spread} hypothesis $|e^\mu|>\left(\rho\left(e^{\D_{x_0}f}\right)\right)^{\infty}\coloneqq 0$ is automatically satisfied, so no assumption is needed on the spectral spread in Prop.~\ref{prop:koopman-cka-fix} (and similarly for Prop.~\ref{prop:koopman-cka-per} below).
	We need only assume that $\mu\neq m\cdot \lambda$ for any $m$ with $\sum_{i}m_i\geq 2$, and this is implied by $\infty$-nonresonance of $(e^\mu, e^{\D_{x_0}f}$) (to see this, take the logarithm of \eqref{eq:nonres-lem}).
	Therefore, Theorem \ref{th:main-thm} and Prop.~\ref{prop:koopman-cka-fix} imply that, for a ``typical'' vector field $f\in \Xfinf$, a unique principal eigenfunction exists for every left eigenvector of $\D_{x_0}f$, and these eigenfunctions are given by a limiting procedure.
	Similar remarks for $f\in \Xpinf$ follow from Theorem \ref{th:main-thm} and Prop.~\ref{prop:koopman-cka-per} below.
\end{rem}

The following is \citet[Prop.~7]{kvalheim2019existence}.
\begin{prop}\label{prop:koopman-cka-per}
	Fix $k\in \N_{\geq 1}\cup \{\infty\}$ and let $\Phi\colon Q \times \R \to Q$ be a $\Cont^k$ flow with $Q$ the basin of an attracting hyperbolic nonstationary $\tau$-periodic orbit with image $\Gamma \subset Q$, where $n+1\coloneqq \dim(Q)\geq 2$.
	Fix $x_0\in \Gamma$ and let $E^s_{x_0}$ be the unique $\D_{x_0}\Phi^\tau$-invariant subspace complementary to $\T_{x_0} \Gamma$.
    Assume the spectral radius $\rho\left(\D_{x_0}\Phi^\tau|_{E^s_{x_0}}\right)<1$ satisfies 
    $$|e^{\mu\tau}|>\left(\rho\left(\D_{x_0}\Phi^\tau|_{E^s_{x_0}}\right)\right)^{k}$$
    in all of the following statements.
    
	\emph{\textbf{Uniqueness of Koopman eigenvalues.}}
    Let $\psi_1\in \Cont^k(Q,\C)$ be any Koopman eigenfunction with eigenvalue $\mu \in \C$.
    Then there exists $m = (m_1,\ldots,m_n)\in \N_{\geq 0}^n$ such that
	$$e^{\mu \tau} =  e^{(m\cdot \lambda)\tau},$$ where  $e^{\lambda_1\tau},\ldots,e^{\lambda_n\tau}$ are the eigenvalues of $\D_{x_0} \Phi^\tau|_{E^s_{x_0}}$ repeated with multiplicities and $\lambda \coloneqq (\lambda_1,\ldots, \lambda_n)$.

	\emph{\textbf{Existence and uniqueness of principal eigenfunctions.}} Assume that $(e^{\mu\tau}, \D_{x_0}\Phi^\tau|_{E^s_{x_0}})$ is $k$-nonresonant. 
	Let $w\colon E^s_{x_0} \to \C$ be a left eigenvector of $\D_{x_0}\Phi^\tau|_{E^s_{x_0}}$ with eigenvalue $e^{\mu \tau}$.
	Then there exists a unique principal eigenfunction $\psi\in \Cont^k(Q,\C)$ for $\Phi$ with eigenvalue $\mu$ satisfying $\D_{x_0}\psi|_{E^s_{x_0}}= w$.
	Additionally, if $\mu$ and $w$ are real, then $\psi\colon Q\to \R \subset \C$ is real.

\end{prop}
\begin{rem}
The uniqueness statements of Prop.~ \ref{prop:koopman-cka-fix} and \ref{prop:koopman-cka-per} are fairly sharp; see \citet[Ex.~2]{kvalheim2019existence}.
\end{rem}

\begin{rem}
	Prop.~\ref{prop:koopman-cka-fix} can be used to guarantee convergence of Laplace averages \citep{mauroy2013isostables}; see \citet[Rem.~14]{kvalheim2019existence}.
	\citet[Rem.~15]{kvalheim2019existence} relates Prop.~\ref{prop:koopman-cka-fix} and \ref{prop:koopman-cka-per} to the literature on isostables and isostable coordinates.
	\citet[Rem.~16]{kvalheim2019existence} relates Prop.~\ref{prop:koopman-cka-fix} to the principal eigenfunctions of \cite{mohr2016koopman}.
\end{rem}

	\subsection{Classification of all $\Cont^\infty$ Koopman eigenfunctions}\label{sec:classify}

	To improve the readability of Theorems \ref{th:classify-point} and \ref{th:classify-per} below, we introduce the following multi-index notation.
	We define an $n$-dimensional multi-index to be an $n$-tuple $i = (i_1,\ldots, i_n)\in \N^n_{\geq 0}$ of nonnegative integers, and define its sum to be $|i|\coloneqq i_1+\cdots +  i_n$.
	For a multi-index $i\in \N^n_{\geq 0}$ and $z = (z_1,\ldots, z_n)\in \C^n$, we define $z^{[i]}\coloneqq z_1^{i_1}\cdots z_n^{i_n}.$
	Given a $\C^n$-valued function $\psi = (\psi_1,\ldots, \psi_n)\colon Q\to \C^n$, we define $\psi^{[i]}\colon Q\to \C$ via $\psi^{[i]}(x)\coloneqq (\psi(x))^{[i]}$ for all $x\in Q$.
	We also define the complex conjugate of $\psi = (\psi_1,\ldots, \psi_n)$ element-wise: $\bar{\psi}\coloneqq (\bar{\psi}_1,\ldots, \bar{\psi}_n).$ 
	The following result follows from \citet[Rem.~3, Ex.~1, Prop.~2, Thm~3]{kvalheim2019existence}.

\begin{thm}[Classification for a point attractor]\label{th:classify-point}
	Let $f$ be a $\Cont^\infty$ vector field on $Q$ with $Q$ the basin of an attracting hyperbolic equilibrium $x_0\in Q$ for the flow of $f$, where $n\coloneqq \dim(Q)\geq 1$.	
	Assume that $\D_{x_0}f$ is diagonalizable over $\C$ with eigenvalues $\lambda\coloneqq (\lambda_1,\ldots, \lambda_n)$ repeated with multiplicities and that $\lambda_j\neq m\cdot \lambda$ for all $j$ and $m\in \N_{\geq 0}$ with $|m|\geq 2$. 
	It follows that there exists an $n$-tuple $\psi = (\psi_1,\ldots,\psi_n)$ of $\Cont^\infty$ principal eigenfunctions such that $\psi\colon Q\to \psi(Q)\subset \C^n$ is a diffeomorphism onto an $\R$-linear subspace of $\C^n$, and every $\Cont^\infty$ Koopman eigenfunction $\varphi$ is a finite linear combination of products of the $\psi_i$ and their complex conjugates $\bar{\psi}_i$:
	\begin{equation}\label{eq:th-classify-expansion}
	\varphi = \sum_{|\ell|+|m| \leq  k}c_{\ell,m}\psi^{[\ell]}\bar{\psi}^{[m]}
	\end{equation}
	for some $k \in \N_{\geq 0}$ and some coefficients $c_{\ell,m}\in \C$.
\end{thm}

\begin{rem}
	Theorem \ref{th:classify-point} goes beyond Prop.~\ref{prop:koopman-cka-fix} by completely classifying \emph{all} $\Cont^\infty$ eigenfunctions rather than just the principal ones.
	On the other hand, Theorem \ref{th:classify-point} requires the stronger hypothesis that $\D_{x_0}\Phi^1$ is diagonalizable over $\C$.
	However, Theorem \ref{th:main-thm} still implies that a ``typical'' vector field in $\Xfinf$ satisfies the hypotheses of Theorem \ref{th:classify-point} (because $\infty$-nonresonance of $(e^{\D_0f},e^{\D_0f})$ implies the condition involving ``$\lambda_j\neq m\cdot \lambda$''; cf. Rem.~\ref{rem:point-cinf-case}).
	Similarly, Theorem \ref{th:classify-per} below yields a complete classification of \emph{all} $\Cont^\infty$ eigenfunctions for ``typical'' vector fields in $\Xpinf$.
\end{rem}

Consider a $\Cont^\infty$ flow with $Q$ the basin of an attracting hyperbolic nonstationary $\tau$-periodic orbit.
As discussed in \cite{kvalheim2019existence}, there exists a $\Cont^\infty$ Koopman eigenfunction $\psi_{\theta}$ with eigenvalue $\mu = i\frac{2\pi }{\tau}$, where $i=\sqrt{-1}$, and this eigenfunction is unique modulo scalar multiplication (cf. \citet{mauroy2012use}).
The following result, which follows from \citet[Prop.~3, Thm~4]{kvalheim2019existence}, involves $\psi_{\theta}$.

\begin{thm}[Classification for a limit cycle attractor]\label{th:classify-per}
	\strut\\
	Let $\Phi\colon Q \times \R \to Q$ be the flow of a $\Cont^\infty$ vector field with $Q$ the basin of an attracting hyperbolic nonstationary $\tau$-periodic orbit with image $\Gamma\subset  Q$, where $n+1\coloneqq \dim(Q)\geq 2$.
	Fix $x_0 \in \Gamma$ and denote by $E^s_{x_0}$ the unique $\tau$-invariant subspace complementary to $\T_{x_0}\Gamma$.
	Assume that $\D_{x_0}\Phi^\tau$ is diagonalizable over $\C$ and that $(\D_{x_0}\Phi^\tau|_{E^s_{x_0}},\D_{x_0}\Phi^\tau|_{E^s_{x_0}})$ is $\infty$-nonresonant.
    It follows that there exists an $n$-tuple $\psi = (\psi_1,\ldots,\psi_n)$ of $\Cont^\infty$ principal eigenfunctions such that $(\psi,\psi_\theta)\colon Q\to \psi(Q)\subset \C^{n+1}$ is a diffeomorphism onto a $\Cont^\infty$ properly embedded submanifold of $\C^{n+1}$, and every $\Cont^\infty$ Koopman eigenfunction $\varphi$ is a finite linear combination of products of $\psi_\theta$ with products of the $\psi_i$ and conjugates $\bar{\psi}_i$:
	\begin{equation}
	\varphi = \sum_{|\ell|+|m|\leq k}c_{\ell,m}\psi^{[\ell]}\bar{\psi}^{[m]} \psi_\theta^{j_{\ell,m}}
	\end{equation}
	for some $k\in \N_{\geq 0}$, coefficients $c_{\ell,m} \in \C$, and $j_{\ell,m} \in \Z$.
\end{thm}

\begin{ack}
We thank George Haller, Hoon Hong, Alexandre Mauroy, Igor Mezi\'{c}, Jeff Moehlis, Ryan Mohr, Corbinian Schlosser, and Dan Wilson for useful discussions related to this work.
\end{ack}

\bibliography{ref}             

\begin{thebibliography}{31}
\providecommand{\natexlab}[1]{#1}
\providecommand{\url}[1]{\texttt{#1}}
\providecommand{\urlprefix}{URL }
\expandafter\ifx\csname urlstyle\endcsname\relax
  \providecommand{\doi}[1]{doi:\discretionary{}{}{}#1}\else
  \providecommand{\doi}{doi:\discretionary{}{}{}\begingroup
  \urlstyle{rm}\Url}\fi

\bibitem[{Blum-Smith and Coskey(2017)}]{blumsmith2017tft}
Blum-Smith, B. and Coskey, S. (2017).
\newblock The fundamental theorem on symmetric polynomials:
  History{\textquotesingle}s first whiff of {G}alois theory.
\newblock \emph{The College Mathematics Journal}, 48(1), 18--29.
\newblock \doi{10.4169/college.math.j.48.1.18}.

\bibitem[{Budi{\v{s}}i{\'c} et~al.(2012)Budi{\v{s}}i{\'c}, Mohr, and
  Mezi{\'c}}]{budivsic2012applied}
Budi{\v{s}}i{\'c}, M., Mohr, R., and Mezi{\'c}, I. (2012).
\newblock Applied {K}oopmanism.
\newblock \emph{Chaos: An Interdisciplinary Journal of Nonlinear Science},
  22(4), 047510.

\bibitem[{Eldering et~al.(2018)Eldering, Kvalheim, and
  Revzen}]{kvalheim2018global}
Eldering, J., Kvalheim, M., and Revzen, S. (2018).
\newblock Global linearization and fiber bundle structure of invariant
  manifolds.
\newblock \emph{Nonlinearity}, 31(9), 4202--4245.

\bibitem[{Hadamard(1902)}]{hadamard1902problemes}
Hadamard, J. (1902).
\newblock Sur les probl{\`e}mes aux d{\'e}riv{\'e}es partielles et leur
  signification physique.
\newblock \emph{Princeton University Bulletin}, 49--52.

\bibitem[{Hirsch(1994)}]{hirsch1976differential}
Hirsch, M.W. (1994).
\newblock \emph{Differential topology}, volume~33 of \emph{Graduate Texts in
  Mathematics}.
\newblock Springer-Verlag, New York.
\newblock Corrected reprint of the 1976 original.

\bibitem[{Khalil(2002)}]{khalil2002nonlinear}
Khalil, H.K. (2002).
\newblock \emph{Nonlinear systems}.
\newblock Prentice Hall.

\bibitem[{{K}oopman(1931)}]{koopman1931hamiltonian}
{K}oopman, B.O. (1931).
\newblock Hamiltonian systems and transformation in {H}ilbert space.
\newblock \emph{Proceedings of the National Academy of Sciences of the United
  States of America}, 17(5), 315.

\bibitem[{Kvalheim and Revzen(2019)}]{kvalheim2019existence}
Kvalheim, M.D. and Revzen, S. (2019).
\newblock Existence and uniqueness of global {K}oopman eigenfunctions for
  stable fixed points and periodic orbits.
\newblock \emph{arXiv preprint arXiv:1911.11996}.

\bibitem[{Lan and Mezi\'c(2013)}]{lan2013linearization}
Lan, Y. and Mezi\'c, I. (2013).
\newblock Linearization in the large of nonlinear systems and {K}oopman
  operator spectrum.
\newblock \emph{Phys. D}, 242, 42--53.
\newblock \doi{10.1016/j.physd.2012.08.017}.
\newblock \urlprefix\url{http://dx.doi.org/10.1016/j.physd.2012.08.017}.

\bibitem[{Lee(2013)}]{lee2013smooth}
Lee, J.M. (2013).
\newblock \emph{Introduction to smooth manifolds}.
\newblock Springer, 2 edition.
\newblock \doi{10.1007/978-1-4419-9982-5}.

\bibitem[{Mauroy and Mezi{\'c}(2012)}]{mauroy2012use}
Mauroy, A. and Mezi{\'c}, I. (2012).
\newblock On the use of {F}ourier averages to compute the global isochrons of
  (quasi) periodic dynamics.
\newblock \emph{Chaos: An Interdisciplinary Journal of Nonlinear Science},
  22(3), 033112.

\bibitem[{Mauroy and Mezi{\'c}(2016)}]{mauroy2016global}
Mauroy, A. and Mezi{\'c}, I. (2016).
\newblock Global stability analysis using the eigenfunctions of the {K}oopman
  operator.
\newblock \emph{IEEE Transactions on Automatic Control}, 61(11), 3356--3369.

\bibitem[{Mauroy et~al.(2013)Mauroy, Mezi{\'c}, and
  Moehlis}]{mauroy2013isostables}
Mauroy, A., Mezi{\'c}, I., and Moehlis, J. (2013).
\newblock Isostables, isochrons, and {K}oopman spectrum for the action--angle
  representation of stable fixed point dynamics.
\newblock \emph{Physica D: Nonlinear Phenomena}, 261, 19--30.

\bibitem[{Mezi\'c(1994)}]{mezic1994geometrical}
Mezi\'c, I. (1994).
\newblock \emph{On the geometrical and statistical properties of dynamical
  systems: theory and applications}.
\newblock Ph.D. thesis, California Institute of Technology.

\bibitem[{Mezi{\'c}(2005)}]{mezic2005spectral}
Mezi{\'c}, I. (2005).
\newblock Spectral properties of dynamical systems, model reduction and
  decompositions.
\newblock \emph{Nonlinear Dynamics}, 41(1-3), 309--325.

\bibitem[{Mezi{\'c}(2019)}]{mezic2019spectrum}
Mezi{\'c}, I. (2019).
\newblock Spectrum of the {K}oopman operator, spectral expansions in functional
  spaces, and state space geometry.
\newblock \emph{arXiv preprint arXiv:1702.07597}.

\bibitem[{Mezi{\'c} and Banaszuk(2004)}]{mezic2004comparison}
Mezi{\'c}, I. and Banaszuk, A. (2004).
\newblock Comparison of systems with complex behavior.
\newblock \emph{Physica D: Nonlinear Phenomena}, 197(1-2), 101--133.

\bibitem[{Mityagin(2015)}]{mityagin2015zero}
Mityagin, B. (2015).
\newblock The zero set of a real analytic function.
\newblock \emph{arXiv preprint arXiv:1512.07276}.

\bibitem[{Mohr and Mezi{\'c}(2016)}]{mohr2016koopman}
Mohr, R. and Mezi{\'c}, I. (2016).
\newblock {K}oopman principle eigenfunctions and linearization of
  diffeomorphisms.
\newblock \emph{arXiv preprint arXiv:1611.01209}.

\bibitem[{Monga et~al.(2019)Monga, Wilson, Matchen, and
  Moehlis}]{monga2019phase}
Monga, B., Wilson, D., Matchen, T., and Moehlis, J. (2019).
\newblock Phase reduction and phase-based optimal control for biological
  systems: a tutorial.
\newblock \emph{Biological Cybernetics}, 113(1-2), 11--46.

\bibitem[{Ott and Yorke(2005)}]{ott2005prevalence}
Ott, W. and Yorke, J. (2005).
\newblock Prevalence.
\newblock \emph{Bulletin of the American Mathematical Society}, 42(3),
  263--290.

\bibitem[{Palis and De~Melo(1982)}]{palis1982geometric}
Palis, J. and De~Melo, W. (1982).
\newblock \emph{Geometric theory of dynamical systems: an introduction}.
\newblock Springer Science \& Business Media.

\bibitem[{Sell(1985)}]{sell1985smooth}
Sell, G.R. (1985).
\newblock Smooth linearization near a fixed point.
\newblock \emph{American Journal of Mathematics}, 1035--1091.

\bibitem[{Shirasaka et~al.(2017)Shirasaka, Kurebayashi, and
  Nakao}]{shirasaka2017phase}
Shirasaka, S., Kurebayashi, W., and Nakao, H. (2017).
\newblock Phase-amplitude reduction of transient dynamics far from attractors
  for limit-cycling systems.
\newblock \emph{Chaos: An Interdisciplinary Journal of Nonlinear Science},
  27(2), 023119.
\newblock \doi{10.1063/1.4977195}.

\bibitem[{Smith and Waltman(1999)}]{smith1999perturbation}
Smith, H.L. and Waltman, P. (1999).
\newblock Perturbation of a globally stable steady state.
\newblock \emph{Proceedings of the American Mathematical Society}, 127(2),
  447--453.

\bibitem[{Sternberg(1957)}]{sternberg1957local}
Sternberg, S. (1957).
\newblock Local contractions and a theorem of {P}oincar{\'e}.
\newblock \emph{American Journal of Mathematics}, 79(4), 809--824.

\bibitem[{Wilson and Ermentrout(2018)}]{wilson2018greater}
Wilson, D. and Ermentrout, B. (2018).
\newblock Greater accuracy and broadened applicability of phase reduction using
  isostable coordinates.
\newblock \emph{Journal of Mathematical Biology}, 76(1-2), 37--66.

\bibitem[{Wilson and Moehlis(2014)}]{wilson2014energy}
Wilson, D. and Moehlis, J. (2014).
\newblock An energy-optimal approach for entrainment of uncertain circadian
  oscillators.
\newblock \emph{Biophysical Journal}, 107(7), 1744--1755.

\bibitem[{Wilson and Moehlis(2016{\natexlab{a}})}]{wilson2016isostable}
Wilson, D. and Moehlis, J. (2016{\natexlab{a}}).
\newblock Isostable reduction of periodic orbits.
\newblock \emph{Physical Review E}, 94(5), 052213.

\bibitem[{Wilson and Moehlis(2016{\natexlab{b}})}]{wilson2016isostable-pde}
Wilson, D. and Moehlis, J. (2016{\natexlab{b}}).
\newblock Isostable reduction with applications to time-dependent partial
  differential equations.
\newblock \emph{Physical Review E}, 94(1), 012211.

\bibitem[{Wilson(1967)}]{wilson1967structure}
Wilson, Jr., F.W. (1967).
\newblock The structure of the level surfaces of a {L}yapunov function.
\newblock \emph{J. Differential Equations}, 3, 323--329.
\newblock \doi{10.1016/0022-0396(67)90035-6}.
\newblock \urlprefix\url{http://dx.doi.org/10.1016/0022-0396(67)90035-6}.

\end{thebibliography}

\appendix

\section{Some background on symmetric polynomials}
\label{app:polynomial}

For the reader's convenience,
this appendix reviews some facts about symmetric polynomials
relevant to Lemma \ref{lem:nonres-full-measure}.
See, e.g., \cite{blumsmith2017tft} for additional background.
The first fact relates a polynomial's coefficients
to certain symmetric polynomials of its roots.
\begin{thm}[Vieta's theorem]
	If $f \in \R[x]$ is a degree $n\geq 1$ monic polynomial
	with roots $\alpha_1,\dots,\alpha_n \in \C$ repeated with multiplicity,
	then
	\begin{align*}
		f(x) &= x^n
			- e_1(\alpha_1,\dots,\alpha_n) x^{n-1}
			+ e_2(\alpha_1,\dots,\alpha_n) x^{n-2}
			\\
			&\quad
			- \cdots
			+ (-1)^n e_n(\alpha_1,\dots,\alpha_n)
		,
	\end{align*}
	where
	\begin{align} \label{eq:elem-sym-poly}
		e_1(\alpha_1,\dots,\alpha_n) &\coloneqq \alpha_1 + \cdots + \alpha_n
		, \\
		e_2(\alpha_1,\dots,\alpha_n) &\coloneqq
			\alpha_1 \alpha_2 + \alpha_1 \alpha_3 + \cdots + \alpha_{n-1} \alpha_n
		, \nonumber \\
		&\;\;\,\vdots \nonumber \\
		e_n(\alpha_1,\dots,\alpha_n) &\coloneqq \alpha_1 \cdots \alpha_n
		, \nonumber
	\end{align}
	are the $n$ \concept{elementary symmetric polynomials}
	in $\alpha_1,\dots,\alpha_n$.
\end{thm}
Since the eigenvalues $\lambda_1,\dots,\lambda_n \in \C$ of $A\in \R^{n\times n}$ are the roots
of the characteristic polynomial $\lambda \mapsto \det(\lambda I_n - A)$,
it follows that its coefficients,
which are polynomials in the matrix entries,
give the elementary symmetric polynomials
in the eigenvalues.
Thus, elementary symmetric polynomials in the eigenvalues
can be obtained directly as polynomials in the matrix entries
\emph{without} computing the eigenvalues.
The fundamental theorem of symmetric polynomials (FTSP)
extends this conclusion to \concept{all symmetric polynomials}
in the eigenvalues,
i.e., to $f \in \R[\lambda_1,\dots,\lambda_n]$ for which
\begin{equation*}
	\forall \sigma \in S_n :
	f(\lambda_1,\dots,\lambda_n)
	\equiv f(\lambda_{\sigma(1)},\dots,\lambda_{\sigma(n)})
	,
\end{equation*}
where $S_n$ is the group of permutations $\sigma$ of $\{1,\ldots,n\}$.
\begin{thm}[FTSP]
	If $f \in \R[\alpha_1,\dots,\alpha_n]$ is symmetric,
	there exists a (unique) polynomial $f_e \in \R[e_1,\dots,e_n]$
	with
	\begin{equation*}
		f(\alpha_1,\dots,\alpha_n)
		\equiv f_e(e_1(\alpha_1,\dots,\alpha_n),\dots,e_n(\alpha_1,\dots,\alpha_n))
		.
	\end{equation*}
\end{thm}
To summarize, any symmetric polynomial in the eigenvalues
is expressible as a polynomial in the matrix entries.

\end{document}